\documentclass[a4paper,12pt]{article}
\usepackage{geometry,latexsym,amssymb,amsmath,amsthm,color}
\usepackage[latin5]{inputenc}
\usepackage{enumerate}
\usepackage{enumitem}
\usepackage[T1]{fontenc}
\usepackage{authblk}
\usepackage[all]{xy}
\usepackage{palatino}
\usepackage{indentfirst}
\usepackage{setspace}
\newtheorem{example}{Example}[section]
\newtheorem{Def}[example]{Definition}
\newtheorem{Exam}[example]{Example}
\newtheorem{Prop}[example]{Proposition}
\newtheorem{Theo}[example]{Theorem}

\newtheorem{Rem}[example]{Remark}
\newtheorem{Cor}[example]{Corollary}
\newenvironment{Prf}{{\bf Proof:} }{\hfill $\Box$
	\mbox{}}
\def\wtilde{\widetilde}
\def\wtilde{\widetilde}
\def\E{\mathsf{E}}

\def\Cov{\mathsf{Cov}}
\def\TC{\mathsf{TC}}
\def\Cat{\mathsf{Cat}}

\def\Ker{\mathsf{Ker}}

\def\St{\mathsf{St}}
\def\C{\mathsf{C}}
\def\XMod{\mathsf{XMod}}

\def\I{\mathsf{I}}

\def\C{\mathsf{C}}
\def\TC{\mathsf{TC}}
\def\E{\mathsf{E}}

\def\qed{\hfill $\Box$}

\begin{document}
\title{Coverings and crossed modules of topological groups with operations}

\author[a]{Osman Mucuk\thanks{E-mail : mucuk@erciyes.edu.tr}}
\author[b]{Tunçar Şahan\thanks{E-mail : tuncarsahan@aksaray.edu.tr}}
\affil[a]{\small{Department of Mathematics, Erciyes University, Kayseri, TURKEY}}
\affil[b]{\small{Department of Mathematics, Aksaray University, Aksaray, TURKEY}}

\date{}

\maketitle
\begin{abstract}
It is a well known result in the covering groups that a subgroup $G$ of the fundamental group at the identity  of a semi-locally simply connected topological group determines a covering morphism of topological groups with characteristic  group $G$. In this paper we generalize this result to a large class  of  algebraic objects called topological groups with operations, including topological groups. We also give  the cover of  crossed modules within  topological groups with operations.
\end{abstract}

\\ {\bf Classification:} 18D35, 22A05, 57M10
	
\section{Introduction}

The theory of covering spaces is one of the most
interesting theories in algebraic topology. It is well known that if $X$ is a topological group, $p\colon \widetilde X\rightarrow X$ is a simply connected covering
map, $e\in X$ is the identity element  and
$\tilde{e}\in \widetilde{X}$ is such that $p(\tilde{e})=e$,
then $\widetilde{X}$ becomes a topological group with identity
$\tilde{e}$ such that $p$ is a morphism of topological groups (see for example \cite{Ch}).

The problem  of universal covers of non-connected topological
groups was first studied  in ~\cite{Ta}. He proved that a
topological group $X$ determines an obstruction class $k_{X}$ in
$H^3(\pi_0(X),\pi_1(X,e))$, and that the vanishing of $k_{X}$ is a
necessary and sufficient condition for the lifting of the group
structure to a universal cover. In ~\cite{Mu1} an analogous
algebraic result was given in terms of crossed modules and group-groupoids, i.e., group objects in the category of groupoids (see also ~\cite{BM1} for
a revised version, which generalizes these results and shows the
relation with the theory of obstructions to extension  for
groups and \cite{Mu-Be-Tu-Na} for the  recently developed  notion of monodromy for topological group-groupoids).

In \cite[Theorem 1]{BS} Brown and Spencer proved that the category of
internal categories within the groups, i.e., group-groupoids,  is equivalent to the
category of crossed modules of groups. Then in \cite[Section 3]{Por}, Porter proved that a similar result holds for certain algebraic categories $\C$, introduced by Orzech \cite{Orz}, which definition was adapted by him and called category of groups with operations. Applying Porter's result, the study of internal category theory in $\C$ was continued in the works of Datuashvili \cite{Wh} and \cite{Kanex}. Moreover, she developed cohomology theory of internal categories, equivalently, crossed modules, in categories of groups with operations \cite{Coh, Cohtr}. In a similar way, the results of  \cite{BS} and  \cite{Por} enabled us to prove that some  properties of covering groups can be generalized to topological groups with operations.

If  $X$ is  a connected topological space which has a universal cover,
$x_0\in X$ and  $G$  is a subgroup of the fundamental group
$\pi_1(X,x_0)$ of $X$ at the point $x_0$, then by \cite[Theorem 10.42]{Rot} we know that there is  a covering map $p\colon
(\widetilde{X}_G,\tilde{x}_0)\rightarrow (X,x_0)$ of pointed spaces, with characteristic group $G$.  In particular if $G$ is singleton, then  $p$  becomes
the universal covering map. Further  if $X$ is a topological group, then $\widetilde{X}_G$ becomes a topological group such that $p$
is a morphism of topological groups. Recently in \cite{Na-Mu1} this method has been applied to topological $R$-modules and obtained a more general result (see also \cite{Na-Mu2} and  \cite{Mu2} for groupoid setting).

The object of this paper is to prove that this result can be generalized to a wide class of algebraic categories, which include categories of  topological groups, topological rings, topological  $R$-modules and alternative topological algebras. This is conveniently handled by working in a category $\TC$.  The method we use is based on that used by Rothman in \cite[Theorem 10.42]{Rot}.
Further we give the cover of  crossed modules within  topological groups with operations.

\section{Preliminaries on groupoids and  covering groups}\label{PrelimCov}

As it is defined in  \cite{Br,Ma} a  groupoid $G$ has a set $G$ of morphisms, which we call just
{\em elements} of $G$, a set $G_0$ of {\em objects} together with maps
$d_0, d_1\colon G\rightarrow G_0$ and $\epsilon \colon G_0 \rightarrow
G$
such that $d_0\epsilon=d_1\epsilon=1_{G_0}$. The
maps
$d_0$, $d_1$ are called {\em initial} and {\em final} point maps respectively
and the map $\epsilon$ is called {\em object inclusion}.
If $a,b\in G$ and $d_1(a)=d_0(b)$, then the {\em composite}
$a\circ b$ exists such that $d_0(a\circ b)=d_0(a)$ and $d_1(a\circ b)=d_1(b)$. So
there exists a partial composition  defined by
$ G_{d_1}\times_{d_0} G\rightarrow G, (a,b)\mapsto a\circ b$, where
$G_{d_1}\times_{d_0} G$ is the pullback of ${d_1}$ and ${d_0}$.
Further, this partial composition is associative, for $x\in G_0$ the element $\epsilon (x)$  acts as the identity, and each element $a$ has an inverse $a^{-1}$ such
that ${d_0}(a^{-1})={d_1}(a)$, ${d_1}(a^{-1})={d_0}(a)$, $a\circ a^{-1}=\epsilon {d_0}(a)$ and
$a^{-1}\circ a=\epsilon {d_1}(a)$. The map $G\rightarrow G$, $a\mapsto a^{-1}$
is called the {\em inversion}.

In a groupoid $G$ for $x,y\in G_0$ we write $G(x,y)$ for the set of all
morphisms with initial point $x$ and final point $y$. According to \cite{Br} for
$x \in G_0$  the {\em star}  of $x$ is defined as $\{a\in G\mid d_0(a)=x\} $ and denoted  as  $\St_G x$.

Let $G$ and $H$ be groupoids. A {\em morphism} from $H$ to $G$ is
a pair of maps  $f\colon H\rightarrow G$ and $f_0\colon
H_0\rightarrow G_0$ such that $d_0 f=f_0 d_0$, $d_1
f=f_0 d_1$, $f\epsilon=\epsilon f_0$ and $f(a\circ b)=f(a)\circ f(b)$ for all $(a,b)\in
H_{d_1}\times_{d_0} H$. For such a morphism we simply write $f\colon
H\rightarrow G$.

Let $p\colon\wtilde G\rightarrow G$ be a morphism of groupoids. Then $p$ is
called a {\em covering morphism} and $\widetilde{G}$  a {\em  covering groupoid} of $G$ if for
each $\tilde x\in {\wtilde G}_0$ the restriction  $\St_{\widetilde{G}} {\tilde{x}}\rightarrow \St_G{p(\tilde x)}$  is  bijective.

We assume the usual theory of covering maps. All spaces $X$ are
assumed to be locally path connected and semi-locally 1-connected,
so that each path component of $X$  admits a simply connected
cover. Recall that a  covering map $p \colon \wtilde{X}
\rightarrow X$ of connected spaces is called  {\it universal}  if
it covers every covering of $X$ in the sense that if $q \colon
\wtilde{Y} \rightarrow X$ is another covering of $X$ then there
exists a map  $r \colon \wtilde{X} \rightarrow \wtilde{Y}$ such
that $p=qr$ (hence $r$ becomes a covering). A covering map $p
\colon \wtilde{X} \rightarrow X$ is called {\it simply connected}
if $\wtilde{X}$ is simply connected. Note that a simply connected
covering is a universal covering.

A  subset $U$ of a space $X$, which has a universal cover, is called  {\em liftable} if it is open, path
connected and lifts to each covering of $X$, that is, if $p \colon
\wtilde{X} \rightarrow X$ is a covering map,  $\imath \colon
U\rightarrow X$ is the inclusion map and $\tilde{x} \in
\wtilde{X}$ such that $p(\tilde x)=x \in U$, then there exists a
map (necessarily unique) $\hat{\imath} \colon U\rightarrow
\wtilde{X}$ such that $p\hat{\imath}=\imath$ and
$\hat{\imath}(x)=\tilde x$. It is an  easy application that  $U$ is liftable if and only if it is open,
path connected and for all $x\in U$, the fundamental group $\pi_{1}(U,x)$ is mapped to the singleton by the morphism
$\imath_{\star}\colon \pi_{1}(U,x)\rightarrow \pi_{1}(X,x) $ induced by the inclusion $\imath\colon (U,x)\rightarrow (X,x)$.

A space $X$ is called {\em semi-locally simply connected} if each point has a liftable neighborhood and {\em locally simply connected} if it has a base of simply connected sets. So a locally simply connected space is also semi-locally simply connected.

For a covering map $p\colon (\widetilde{X},\tilde{x}_0)\rightarrow
(X,x_0)$ of pointed topological spaces, the subgroup
$p_{\star}(\pi_1(\widetilde{X},\tilde{x}_0))$ of
$\pi_1(X,x_0)$ is called {\em characteristic group} of $p$, where
$p_{\star}$ is the morphism induced by $p$ (see  \cite[p.379]{Br} for  the characteristic group of a covering map  in terms of covering morphism of  groupoids ). Two covering maps $p\colon
(\widetilde{X},\tilde{x}_0)\rightarrow (X,x_0)$ and
$q\colon (\widetilde{Y},\tilde{y}_0)\rightarrow (X,x_0)$
are called {\em equivalent} if their characteristic groups are isomorphic,  equivalently  there is a homeomorphism $f\colon
(\widetilde{X},\tilde{x}_0)\rightarrow
(\widetilde{Y},\tilde{y}_0)$ such that $qf=p$.

We recall a  construction from \cite[p.295]{Rot} as follows: Let $X$ be a topological space with a base point $x_0$ and  $G$ a subgroup of $\pi_1(X,x_0)$. Let   $P(X,x_0)$ be  the
set of all paths of $\alpha$ in $X$ with initial point $x_0$. Then the relation defined   on $P(X,x_0)$ by $\alpha\simeq\beta$ if and only if
$\alpha(1)=\beta(1)$ and $[\alpha\circ \beta^{-1}]\in G$,  is an equivalence relation.  Denote the
equivalence relation of $\alpha$ by $\langle \alpha\rangle _G$ and
define $\widetilde{X}_G$ as the set of all such equivalence
classes of the paths in  $X$ with initial point $x_0$. Define a function $p\colon \widetilde{X}_G\rightarrow
X$ by $p(\langle \alpha\rangle _G)=\alpha(1)$. Let $\alpha_0$ be
the constant path at $x_0$ and $\tilde{x}_0=\langle
\alpha_0\rangle _G\in \widetilde{X}_G$. If  $\alpha\in P(X,x_0)$ and $U$ is  an open
neighbourhood of $\alpha(1)$, then  a path  of  the form
$\alpha \circ \lambda$, where $\lambda$ is a path in $U$ with
$\lambda(0)=\alpha(1)$, is called a {\em continuation} of
$\alpha$. For an $\langle \alpha\rangle _G\in \widetilde{X}_G$ and an open
neighbourhood $U$ of $\alpha(1)$, let $(\langle \alpha\rangle _G,U)=\{\langle\alpha\circ\lambda \rangle_G\colon ~\lambda(I)\subseteq
U\}$. Then the subsets
$(\langle \alpha\rangle _G, U)$ form a basis for a topology on
$\widetilde{X}_G$ such that the map $p\colon
(\widetilde{X}_G,\tilde{x}_0)\rightarrow (X,x_0)$ is
continuous.

In Theorem \ref{Lifttopmodul} we generalize the following result to topological groups with operations.
\begin{Theo}\label{CovCorrG}{\em \cite[Theorem 10.34]{Rot}} Let $(X,x_0)$ be a pointed topological space and
	$G$ a subgroup of $\pi_1(X,x_0)$. If $X$ is connected, locally
	path connected and semi-locally simply connected, then $p\colon
	(\widetilde{X}_G,\tilde{x}_0)\rightarrow (X,x_0)$ is a
	covering map with characteristic  group $G$ .\end{Theo}

\begin{Rem}\label{Rem} {\em Let $X$ be a connected, locally path connected
		and semi locally simply connected topological space and $q\colon
		(\widetilde{X},\tilde{x}_0)\rightarrow (X,x_0)$  a
		covering map. Let $G$ be the  characteristic group  of $q$. Then
		the covering map $q$ is equivalent to the covering map $p\colon
		(\widetilde{X}_G,\tilde{x}_0)\rightarrow (X,x_0)$
		corresponding to $G$. }\end{Rem} So  from Theorem \ref{CovCorrG}
the following result is obtained.
\begin{Theo}\label{TheoRotman} {\em \cite[Theorem 10.42]{Rot}} Suppose that  $X$ is  a connected,
	locally path connected and semi-locally simply connected
	topological group. Let  $e\in X$ be the identity element and
	$p\colon (\widetilde{X},\tilde{e})\rightarrow (X,e)$ a
	covering map. Then the group structure of $X$ lifts to
	$\widetilde{X}$, i.e.,  $\widetilde{X}$ becomes a topological group
	such that $\tilde{e}$ is identity  and $p\colon
	(\widetilde{X},\tilde{e})\rightarrow (X,e)$ is a morphism
	of topological groups.\end{Theo}

\section{Universal covers of topological groups with operations}

In this section we apply the methods of Section \ref{PrelimCov} to the topological groups with operations and obtain parallel results.

The idea  of the definition of categories of groups with
operations comes from  \cite{Hig} and  \cite {Orz} (see also \cite{Orz2})
and the definition below is from \cite{Por} and  \cite[p.21]{Tamar}, which is adapted from \cite {Orz}.

\begin{Def} \label{Defgroupwithoperation}{\em
		Let $\C$ be a  category of groups with a set  of  operations $\Omega$ and with a set $\E$  of identities such that $\E$ includes the group laws, and the following conditions hold for the set $\Omega_i$  of $i$-ary operations in $\Omega$:
		
		(a) $\Omega=\Omega_0\cup\Omega_1\cup\Omega_2$;
		
		(b) The group operations written additively $0,-$ and $+$ are
		the  elements of $\Omega_0$, $\Omega_1$ and
		$\Omega_2$ respectively. Let $\Omega_2'=\Omega_2\backslash \{+\}$,
		$\Omega_1'=\Omega_1\backslash \{-\}$ and assume that if $\star\in
		\Omega_2'$, then $\star^{\circ}$ defined by
		$a\star^{\circ}b=b\star a$ is also in $\Omega_2'$. Also assume
		that $\Omega_0=\{0\}$;
		
		(c) For each   $\star \in \Omega_2'$, $\E$ includes the identity
		$a\star (b+c)=a\star b+a\star c$;
		
		(d) For each  $\omega\in \Omega_1'$ and $\star\in \Omega_2' $, $\E$
		includes the identities  $\omega(a+b)=\omega(a)+\omega(b)$ and
		$\omega(a)\star b=\omega(a\star b)$.
		
		The category $\C$ satisfying the conditions (a)-(d) is called a {\em category of groups with operations}}.\end{Def}

In the paper from now on  $\C$ will denote the category of groups with operations.

A {\em  morphism} between any two objects of $\C$ is a group homomorphism, which preserves the operations of $\Omega_1'$ and $\Omega_2'$.

\begin{Rem}{\em  The set $\Omega_0$ contains exactly one element, the
		group identity; hence for instance the category of associative rings with unit is not a category of  groups with operations.}\end{Rem}

\begin{Exam}{\em The categories of  groups, rings generally  without identity, $R$-modules,  associative, associative commutative, Lie, Leibniz,
		alternative algebras are examples of categories of   groups with operations.\qed}
\end{Exam}

\begin{Rem}{\em  The set $\Omega_0$ contains exactly one element, the
		group identity; hence for instance associative rings with unit are
		not groups with operations}\end{Rem}

The  category of topological groups with operations is defined in \cite{Ak-Na-Mu-Tu} as follows:

\begin{Def}{\em  A category $\TC$  of topological groups with a set $\Omega$ of continuous operations  and with a set $\E$ of identities such  that $\E$ includes the group laws such that the conditions (a)-(d) in Definition \ref{Defgroupwithoperation} are satisfied, is called a {\em   category of topological groups with operations} and the object of  $\TC$ are called {\em topological groups with operations}}.\end{Def}

In the rest of the paper $\TC$ will denote a category of topological groups with operations.

A {\em  morphism} between any two objects of $\TC$ is a continuous group homomorphism, which preserves the operations in $\Omega_1'$ and $\Omega_2'$.

The categories of topological groups, topological  rings,  topological $R$-modules and alternative topological  algebras are examples of categories of  topological groups with operations.

\begin{Prop}  If  $X$ is a topological group with operations, then the fundamental group $\pi_1(X,0)$ becomes a group with operations.\end{Prop}
\begin{Prf} Let $X$ be an object of $\TC$ and  $P(X,0)$  the
	set of all paths in $X$ with initial point $0$ as described in Section 1.
	There are binary operations on $P(X,0)$ defined by
	\begin{align}\label{1}
		(\alpha \star \beta)(t)=\alpha(t)\star \beta(t)
	\end{align}
	for $\star\in \Omega_2$ and $t\in \I$, unit interval, and unary operations defined by
	\begin{align}\label{2}
		(\omega \alpha)(t)=\omega(\alpha(t))
	\end{align}
	for $\omega\in \Omega_1$.   Hence  the operations $(\ref{1})$ induce binary operations on $\pi_1(X,0)$ defined by
	\begin{align}\label{defeqrel}
		[\alpha] \star [\beta]=[\alpha\star \beta]
	\end{align}
	for  $[\alpha],[\beta]\in  \pi_1(X,0)$. Since the binary operations $\star$  in $\Omega_2$ are continuous it follows that the binary operations $(\ref{defeqrel})$ are well defined. Similarly the operations  $(\ref{2})$ reduce the unary operations defined by
	\begin{align}\label{4}
		\omega[\alpha] =[\omega\alpha].
	\end{align}
	By the continuity of the unary  operations $\omega\in \Omega_1$, the operations $(\ref{4})$ are also well defined.
	The other details can be checked and so $\pi_1(X,x_0)$  becomes  a group with operations, i.e., an object of $\C$. \end{Prf}

We now generalize Theorem \ref{CovCorrG} to  topological groups with operations.  We first make the following preparation:

Let $X$ be a topological group with operations. By the evaluation  of the compositions  and operations of the paths in  $X$  such that  $\alpha_1(1)=\beta_1(0)$ and $\alpha_2(1)=\beta_2(0)$   at $t\in \I$ ,   we have the following {\em interchange law}  \begin{align}\label{interchange}
	(\alpha_1\circ \beta_1)\star(\alpha_2\circ \beta_2)=(\alpha_1\star \alpha_2)\circ (\beta_1\star \beta_2)
\end{align}
for $\star\in \Omega_2$, where $\circ$ denotes the composition of paths, and
\begin{align}\label{6}
	(\alpha\star \beta)^{-1}=\alpha^{-1}\star\beta^{-1}
\end{align}
for $\alpha,\beta\in P(X,0)$ where, say  $\alpha^{-1}$ is the inverse path  defined by $\alpha^{-1}(t)=\alpha(1-t)$ for $t\in I$. Further we have that
\begin{align}\label{7}
	(\omega \alpha)^{-1}=\omega \alpha^{-1}
\end{align}
\begin{align}\label{8}
	\omega (\alpha\circ \beta)=(\omega \alpha)\circ (\omega\beta)
\end{align}
when $\alpha(1)=\beta(0)$.

Parallel to Theorem \ref{CovCorrG}, in the following theorem we prove a general result for topological groups with operations.
\begin{Theo} \label{Lifttopmodul} Let $X$ be a topological group with operations, i.e., an object of $\TC$ and let  $G$ be the subobject of  $\pi_1(X,0)$. Suppose that the underlying space of $X$ is connected, locally path connected and semi-locally simply connected. Let
	$p\colon (\widetilde{X}_G,\tilde{0})\rightarrow (X,0)$ be
	the covering map corresponding  to  $G$ as a subgroup of the additive group of $\pi_1(X,0)$  by Theorem
	\ref{CovCorrG}. Then the group operations  of $X$ lift to
	$\widetilde{X}_G$, i.e., $\widetilde{X}_G$ is a topological group with operations
	and $p\colon \widetilde{X}_G\rightarrow X$ is a
	morphism of $\TC$.\end{Theo}
\begin{Prf} By the construction of $\widetilde{X}_G$ in Section 1,  $\widetilde{X}_G$ is the set of equivalence classes defined via $G$.   The  binary  operations  on $P(X,0)$  defined by $(\ref{1})$  induce binary operations
	\begin{align}\label{9}
		\langle \alpha\rangle _G\star\langle \beta\rangle _G=\langle
		\alpha\star \beta\rangle _G
	\end{align}
	and the unary operations on  $P(X,x_0)$  defined by $(\ref{2})$  induce unary  operations
	\begin{align}\label{10}
		\omega\langle\alpha\rangle_G=\langle\omega\alpha\rangle_G \end{align} on $\widetilde{X}_G$.
	
	We now prove that these operations $(\ref{9})$ and $(\ref{10})$ are well defined: For   $\star\in \Omega_2$ and  the paths $\alpha, \beta, \alpha_1, \beta_1\in P(X,0)$ with $\alpha(1)=\alpha_1(1)$ and $\beta(1)=\beta_1(1)$, we have that
	\begin{align*}
		[(\alpha\star \beta)\circ (\alpha _{1}\star \beta_{1})^{-1}]&=[(\alpha\star \beta)\circ (\alpha_1^{-1}\star \beta_1^{-1})] \tag{by \ref{6}}\\
		&=[(\alpha\circ \alpha_1^{-1})\star(\beta\circ \beta_1^{-1})] \tag{by \ref{interchange}}\\
		&=[\alpha\circ \alpha_1^{-1}]\star [\beta\circ \beta_1^{-1}] \tag{by  \ref{defeqrel}}
	\end{align*}
	So if  $\alpha_{1}\in \langle \alpha\rangle _G $ and $\beta_{1}\in \langle \beta\rangle _G $ , then $[\alpha\circ \alpha_1^{-1}]\in G$ and
	$[\beta\circ \beta_1^{-1}]\in G$.  Since $G$ is a subobject of $\pi_1(X,0)$, we have that $[\alpha\circ \alpha_1^{-1}]\star [\beta\circ \beta_1^{-1}]\in G$. Therefore the binary operations $(\ref{9})$ are well defined.
	
	Similarly for the paths $\alpha,\alpha_1 \in P(X,0)$ with $\alpha(1)=\alpha_1(1)$  and $\omega\in \Omega_1$ we have that
	\begin{align*}
		[( \omega \alpha) \circ  (\omega \alpha_1)^{-1})]&= [( \omega \alpha) \circ  (\omega \alpha_1^{-1})] \tag{by  \ref{7}}\\
		&=[( \omega (\alpha \circ\alpha_1^{-1})] \tag{by \ref{8}}\\
		&=\omega [\alpha \circ\alpha_1^{-1}] \tag{by \ref{4}}
	\end{align*}
	Since $G$ is a subobject of $\pi_1(X,0)$, if $[\alpha\circ \alpha_1^{-1}]\in G$ and $\omega\in \Omega_1$ then $\omega [\alpha\circ \alpha_1^{-1}]\in G$. Hence the unary operations $(\ref{10})$ are also  well defined.
	
	The axioms (a)-(d) of Definition  \ref{Defgroupwithoperation}  for $\widetilde{X}_G$  are satisfied and therefore  $\widetilde{X}_G$ becomes a group with operations. Further by Theorem \ref{CovCorrG},  $p\colon
	(\widetilde{X}_G,\tilde{0})\rightarrow (X,0)$  is a  covering map, $\widetilde{X}_G$ is a topological group and $p$ is a morphism of topological groups. In addition to this we need to prove that   $\widetilde{X}_G$ is  an object of $\TC$ and $p$ is a morphism of $\TC$. To prove that the operations $(\ref{9})$ for $\star\in \Omega_2'$  are continuous let  $\langle
	\alpha\rangle _G,\langle\beta\rangle _G\in \widetilde{X}_G $  and $(W, \langle
	\alpha\star \beta\rangle _G)$  be a basic open neighbourhood of $\langle
	\alpha\star\beta\rangle_G $. Here  $W$ is an open neighbourhood of
	$(\alpha\star\beta)(1)=\alpha(1)\star \beta(1)$. Since the operations  $\star\colon X\times
	X\rightarrow X$ are continuous there are open neighbourhoods $U$
	and $V$ of $\alpha(1)$ and $\beta(1)$ respectively in $X$ such
	that $U\star V\subseteq W$. Therefore
	$(U,\langle \alpha\rangle_G)$ and $ (V,\langle \beta\rangle_G)$ are respectively base open neighbourhoods of $\langle\alpha\rangle_G$ and    $\langle\beta\rangle_G$, and  \[(U,\langle \alpha\rangle_G)\star (V,\langle \beta\rangle_G)\subseteq (W,\langle
	\alpha\star\beta\rangle _G).\] Therefore the binary operations $(\ref{9})$ are continuous
	
	We now prove that the unary operations $(\ref{10})$  for $\omega\in \Omega_1'$ are continuous. For if $(V,\langle \omega\alpha\rangle)$ is a base open neighbourhood of $\langle \omega\alpha\rangle$, then $V$ is an open neighbourhood of $\omega\alpha(1)$ and  since the unary operations  $\omega\colon X\rightarrow X$ are continuous there is an open neighbourhood $U$ of $\alpha(1)$ such that $\omega(U)\subseteq V$. Therefore $(U,\langle \alpha\rangle)$ is an open neighbourhood of $\langle \alpha\rangle$ and $\omega (U,\langle \alpha\rangle)\subseteq (V,\langle\omega \alpha\rangle)$.
	
	Moreover the map  $p\colon \widetilde{X}_G\rightarrow X$   defined by  $p(\langle \alpha \rangle_G)=\alpha(1)$ preserves the operations of $\Omega_2$ and $\Omega_1$. \end{Prf}

From Theorem \ref{Lifttopmodul} the following result can be restated.
\begin{Theo}\label{TheoLiftingTopgpwithopretin} Suppose that  $X$ is a topological group with operations whose underlying
	space is connected, locally path connected and semi-locally simply
	connected. Let  $p\colon
	(\widetilde{X},\tilde{0})\rightarrow (X,0)$ be a covering map
	such that $\widetilde{X}$ is path connected and the characteristic
	group $G$ of $p$ is a subobject of   $\pi_1(X,0)$. Then the group operations of $X$ lifts to $\widetilde{X}$.\end{Theo}
\begin{Prf}   By assumption the characteristic group $G$ of the covering map $p\colon
	(\widetilde{X},\tilde{0})\rightarrow (X,0)$  is a subobject of $\pi_1(X,0)$. So by Remark
	\ref{Rem},  we can assume that $\widetilde{X}=\widetilde{X}_G$ and hence  by Theorem \ref{Lifttopmodul}, the group operations of $X$ lift to $\widetilde{X}$ as required.
\end{Prf}

In particular, in  Theorem \ref{Lifttopmodul} if the subobject $G$ of $\pi_1(X,0)$ is chosen to be
the singleton, then   the following corollary  is obtained.
\begin{Cor} Let $X$ be a topological group with operations such that the underlying space of $X$ is connected,
	locally path connected and semi-locally simply connected. Let  $p\colon
	(\widetilde{X},\tilde{0})\rightarrow (X,0)$ be  a
	universal covering map.  Then the group structures  of $X$ lifts to $\widetilde{X}$.\end{Cor}

The following proposition is useful for Theorem  \ref{Liftofinclusion}.
\begin{Prop}\label{Liftablenbd}
	Let $X$ be  a topological group with operations and $V$ a liftable neighbourhood of $0$ in $X$. Then there is a liftable neighbourhood
	$U$ of $0$ in $X$  such that $U\star U\subseteq V$  for $\star\in
	\Omega_2 $. \end{Prop}
\begin{Prf}
	Since $X$ is a topological group  with operations and hence the binary operations  $\star\in \Omega_2$ are  continuous, there is an open neighbourhood  $U$ of
	$0$ in $X$ such that $U\star U\subseteq V$.  Further if $V$ is
	liftable, then $U$ can be chosen as liftable. For if $V$ is
	liftable, then for each $x\in U$, the fundamental group
	$\pi_{1}(U,x)$ is mapped to the singleton by the morphism induced
	by the inclusion map $\imath \colon U\rightarrow X$. Here $U$ is not
	necessarily path connected and hence not necessarily liftable. But
	since the path component $C_{0}(U)$ of $0$ in $U$ is liftable and
	satisfies these conditions,  $U$ can be replaced by  the the path component $C_{0}(U)$ of $0$ in $U$ and assumed that  $U$ is liftable.
\end{Prf}

\begin{Def}{\em
		Let $X$ and $Y$ be topological groups with operations and $U$ an open neighbourhood of  $0$ in $X$. A
		continuous map $\phi \colon U\rightarrow S$ is called a {\em local
			morphism} in $\TC$ if $\phi(a\star b)=\phi(a)\star \phi(b)$ when  $a,b \in
		U$ such that $a\star b\in U$ for $\star \in \Omega_2$.}
\end{Def}
\begin{Theo}\label{Liftofinclusion}
	Let $X$ and $\widetilde{X}$ be  topological groups with operations and $q \colon \widetilde{X}\rightarrow X$ a morphism of $\TC$, which is a covering map.
	Let $U$ be an open, path connected neighbourhood of  $0$ in $X$ such that
	for each $\star\in \Omega_2$,  the set  $U\star U$ is  contained
	in a liftable neighbourhood $V$ of $0$ in $X$. Then the inclusion
	map $i\colon U\rightarrow X$ lifts to a local morphism
	$\hat{\imath} \colon U\rightarrow \widetilde{X}$ in $\TC$.
\end{Theo}
\begin{Prf}
	Since $V$ lifts to $\widetilde{X}$, then $U$ lifts to
	$\widetilde{X}$ by $\hat{\imath} \colon U\rightarrow
	\widetilde{X}$. We now prove that $\hat{\imath}$ is a local
	morphism of topological groups with operations. We know by the
	lifting theorem that $\hat{\imath} \colon U\rightarrow
	\widetilde{X}$ is continuous. Let $a,b\in U$ be such that for each
	$\star\in \Omega_2$,  $a\star b\in U$. Let $\alpha$ and $\beta$ be the
	paths from $0$ to $a$ and $b$ respectively in $U$. Let $\gamma=\alpha\star
	\beta$. So  $\gamma$ is a  path from  $0$ to $a\star b$.  Since $U\star U
	\subseteq V$,  the paths $\gamma$ is  in V.  So the paths $\alpha,\beta$ and $\gamma$
	lift to $\widetilde{X}$. Suppose that $\tilde{\alpha}$, $\tilde{\beta}$ and  $\tilde{\gamma}$  are the
	liftings of $\alpha$, $\beta$ and  $\gamma$  in $\widetilde{X}$ respectively.
	Then we have
	\[       q(\wtilde{\gamma})=\gamma=\alpha\star\beta= q(\wtilde{\alpha})\star q(\wtilde{\beta}).\]
	But  $q$ is a morphism of topological group with operations and so
	we have,
	\[     q(\wtilde{\alpha}\star \wtilde{\beta})=q(\wtilde{\alpha})\star q(\wtilde{\beta})\] for $\star\in \Omega_2$.
	Since the paths  $\tilde{ \gamma}$ and
	$\tilde{\alpha}\star\tilde{\beta}$ have the initial point
	$\tilde{0}\in \widetilde{X}$,  by the unique path lifting
	\[\wtilde \gamma=\wtilde{\alpha}\star\tilde{\beta} \]
	On evaluating these paths at $1\in \I$ we
	have
	\[     \hat{ {\imath}}(a\star b) =\hat{{\imath}}(a)\star\hat{{\imath}}(b). \]
\end{Prf}

\section{ Covers of crossed modules within topological groups with operations}

If $A$ and $B$ are objects of $\C$, an {\em extension} of $A$ by $B$ is
an exact sequence
\begin{align}
	0\longrightarrow A\stackrel{\imath}\longrightarrow E\stackrel{p}\longrightarrow B\longrightarrow 0
\end{align}
in which $p$ is surjective and $\imath$ is the kernel of $p$.  It is {\em split } if there is a morphism $s\colon  B \to E$ such
that $p s = \imath d _B$.  A split extension of $B$ by $A$ is called a {\em  $B$-structure} on $A$.  Given  such a $B$-structure on $A$ we  get actions of $B$ on $A$ corresponding to the operations in $\C$. For any $b\in B$, $a\in A$ and $\star\in \Omega'_2$ we have the actions called {\em derived actions} by Orzech \cite[p.293]{Orz}
\begin{equation} \label{eq12} \begin{array}{rcl}
		b\cdot a & = &s(b)+a-s(b)\\
		b\star a  & = &s(b)\star a.
	\end{array}\end{equation}
	In addition to this  we note that topologically if an exact sequence $(11)$ in $\TC$ is an split extension, then the derived actions (12) are continuous. So we can state Theorem \cite[Theorem 2.4]{Orz}  in topological case, which is useful for the proof of Theorem \ref{Theocatequivalencess}, as follows.

	\begin{Theo} \label{TeoTopDerivedaction}  A set of actions (one for each operation in $\Omega_2$) is a set of
		continuous derived actions if and only if the semidirect product $B \ltimes
		A$ with underlying set $B \times A$ and operations
		\begin{eqnarray*}
			(b, a)+(b', a') &=& (b + b', a+(b\cdot a')) \\
			(b, a) \star (b', a') &=& (b \star b', a \star a'+ b \star a' +a \star b')
		\end{eqnarray*}
		is an object in $\TC$.\end{Theo}

	The internal category in $\C$ is defined in \cite{Por} as follows. We  follow the notations of Section 1 for groupoids.
	\begin{Def}\label{Internalgpd} {\em  An internal category $C$  in $\C$  is a category in which the initial and final point maps $d_0,d_1\colon C\rightrightarrows
			C_0$, the object inclusion map  $\epsilon\colon C_0\rightarrow C$
			and the partial composition $\circ\colon C_{d_1}\times_{d_0}
			C\rightarrow C,(a,b)\mapsto a\circ b$ are the morphisms in the category
			$\C$.\qed }\end{Def}
	
	Note that since $\epsilon$ is a morphism in $\C$,
	$\epsilon(0)=0$ and that the operation $\circ$ being a morphism
	implies that for all $a,b,c,d\in C$ and $\star\in \Omega_2$,
	\begin{align}\label{13}
		(a\star b)\circ(c\star d)=(a\circ c)\star(b\circ d)
	\end{align}
	whenever one side makes sense. This is called the {\em interchange law} \cite{Por} .
	
	We also note from \cite{Por} that  any internal category in
	$\C$ is an internal groupoid since given $a\in C$,
	$a^{-1}=\epsilon d_1(a)-a+\epsilon d_0(a)$ satisfies $a^{-1}\circ a=\epsilon d_1(a)$ and
	$a\circ a^{-1}=\epsilon d_0(a)$. So we use the term {\em internal groupoid}  rather than  internal category and write  $G$ for an internal groupoid.  For the category of internal groupoids in $\C$ we use   the same notation $\Cat(\C)$ as in \cite{Por}.  Here a {\em morphism}  $f\colon H\rightarrow G$ in $\Cat(\C)$  is morphism of underlying groupoids and a morphism in $\C$.

	In particular if $\C$ is the  category of groups, then an internal groupoid  $G$ in $\C$
	becomes a  group-groupoid and in the case where $\C$ is the category of rings, an internal groupoid in $\C$ is a ring object in the category of groupoids \cite{Mu2}.

	\begin{Def}{\em  An internal groupoid  in the category  $\TC$ of topological groups with operations is called a {\em  topological internal groupoid}.}\end{Def}
	
	So  a topological internal groupoid is a topological groupoid $G$ in which the set of morphisms  and the set $G_0$  of objects are objects of $\TC$ and all structural maps of $G$, i.e, the source and target maps $s,t\colon G\rightarrow G_0$, the object inclusion map  $\epsilon\colon G_0\rightarrow G$ and  the composition map $\circ\colon G_{t}\times_{s} G\rightarrow G$, are  morphisms of  $\TC$.
	
	If $\TC$ is the category of topological groups, then an internal topological groupoid becomes a topological group-groupoid.
	
	For the category of internal topological groupoids in $\TC$ we use   the notation $\Cat(\TC)$.  Here a {\em morphism}  $f\colon H\rightarrow G$ in $\Cat(\TC)$  is morphism of underlying groupoids and a morphism in $\TC$.

	\begin{Theo}\label{Propintgroup} Let    $X$ be an object of $\TC$ such that the underlying space   is locally path connected and  semi-locally  simply connected. Then the fundamental groupoid $\pi X$  is   a topological internal groupoid.
	\end{Theo}
	\begin{Prf} Let $X$ be a topological group with operations as assumed.  By  \cite[Theorem 1]{Br-Da},  $\pi X$ has a topology such that  it is a topological  groupoid.   We know by   \cite[Proposition 3]{Br-Da} that when $X$ and $Y$ are endowed with such topologies, for a continuous map $f\colon X\rightarrow Y$, the induced morphism $\pi(f)\colon \pi X\rightarrow \pi Y$ is also continuous.  Hence the continuous binary operations $\star\colon X\times X \rightarrow X$ for $\star\in \Omega_2$ and  the unary operations $\omega\colon X\rightarrow X$ for $\omega\in \Omega_1$ respectively  induce  continuous binary  operations   $\widetilde{\star}\colon \pi X\times \pi X\rightarrow \pi X$ and unary operations  $\widetilde{\omega}\colon \pi X\rightarrow \pi X$. So the set of morphisms becomes   a topological group with operations.  The groupoid structural maps are morphisms of groups with operations, i.e., preserve the operations. Therefore $\pi X$ becomes a  topological internal groupoid.\end{Prf}

	As similar to the crossed module in $\C$ formulated in \cite[Proposition 2]{Por}, we define a crossed module in $\TC$  as follows:
	
	\begin{Def} \label{Defcrosmod} {\em  A {\em crossed module} in $\TC$ is a morphism  $\alpha\colon A\rightarrow B$ in $\TC$, where $B$ acts topologically  on $A$ (i.e. we have a continuous derived action in $\TC$)  with the conditions for any $b \in B$, $a, a'\in A$, and $\star\in\Omega_2'$:
			\begin{enumerate}
				\item [CM1]$\alpha(b \cdot a) = b + \alpha(a) - b$;
				\item [CM2]$\alpha(a)\cdot a'=a+a'-a$;
				\item [CM3] $\alpha(a)\star a'=a\star a'$;
				\item [CM4] $\alpha(b\star a)=b\star\alpha(a)$ and  $\alpha(a\star b)=\alpha(a)\star b$.
			\end{enumerate}\qed}\end{Def}
	
	A {\em morphism} from  $\alpha\colon A\rightarrow B$ to $\alpha'\colon A'\rightarrow B'$   is a pair
	$f_1\colon A\rightarrow A'$ and  $f_2\colon B\rightarrow B'$ of the morphisms in $\TC$ such that
	\begin{enumerate}[label={ \arabic{*}.} , leftmargin=1cm]
		\item $f_2 \alpha(a)=\alpha' f_1(a)$,
		\item $f_1(b\cdot a)=f_2(b)\cdot f_1(a)$,
		\item $f_1(b\star a)=f_2(b)\star f_1(a)$
	\end{enumerate}
	for any $x\in B$, $a\in A$ and $\star\in\Omega_2'$.  So  we have a category $\XMod(\TC)$ of crossed modules in $\TC$.

	The algebraic case of the following theorem was proved in $\C$ in \cite[Theorem 1]{Por}. We can state the topological version  as follows.
	\begin{Theo} \label{Theocatequivalencess}  The category $\XMod(\TC)$ of crossed modules in $\TC$  and the category $\Cat(\TC)$  of internal groupoids in $\TC$ are equivalent.\end{Theo}
	\begin{Prf} We give a sketch proof  based on  that of algebraic case. A functor $\delta\colon \Cat(\TC)\rightarrow \XMod(\TC)$ is defined as follows:
		For a topological internal groupoid $G$,  let  $\delta(G)$ be the topological crossed module
		$(A,B,d_1)$ in $\TC$,  where $A=\Ker d_0$, $B=G_0$ and $d_1\colon A\rightarrow B$ is the
		restriction of the target point map. Here $A$ and $B$ inherit the structures of topological group with
		operations  from that of $G$, and the target point map
		$d_1\colon A\rightarrow B$ is a morphism in $\TC$. Further the  actions $B\times A\rightarrow A$ on
		the topological group with operations $A$ given by
		\begin{eqnarray*}
			b\cdot a &=& \epsilon(b)+a-\epsilon(b) \\
			b\star a &=& \epsilon(b)\star a
		\end{eqnarray*}
		for  $a\in A$, $b\in B$ are continuous by the continuities  of  $\epsilon$  and the operations in $\Omega_2$; and  the axioms of Definition \ref{Defcrosmod} are satisfied. Thus $(A,B,d_1)$ becomes a crossed module in $\TC$.
		
		Conversely define a functor $\eta\colon \XMod(\TC)\rightarrow \Cat(\TC)$ in the following way. For a crossed module  $(A,B,\alpha)$ in $\TC$,  define a topological internal groupoid   $\eta(A,B,\alpha)$  whose  set of objects is the topological group with operations $B$ and  set of morphisms is the semi-direct product
		$B\ltimes A$ which  is a topological group with operations by Theorem \ref{TeoTopDerivedaction}.
		The  source and target  point maps are defined to be $d_0(b,a)= b$ and  $d_1(b,a)= \alpha(a)+b$  while the object inclusion map and groupoid
		composition is given by $\epsilon(b)=(b,0)$ and \[(b,a)\circ (b_1,a_1)=(b,a_1+a)\] whenever
		$b_1=\alpha(a) b$. These structural maps  are all continuous and therefore $\eta(A,B,\alpha)$ is a topological internal groupoids.
		
		The other details of the proof  is obtained from that of \cite[Theorem 1]{Por}.
		
	\end{Prf}
	
	By Theorem \ref{Theocatequivalencess}, we obtain the cover of a crossed module in $\TC$.  If $f\colon H\rightarrow G$ is a covering morphism in $\Cat(\TC)$ and $(f_1,f_2)$ is the morphism of crossed modules  corresponding to $f$, then $f_1\colon A\rightarrow A'$ is an isomorphism in $\TC$, where $A=\St_H 0$, $A'=\St_G 0$ and $f_1$ is the restriction of $f$.  Therefore we call a morphism $(f_1,f_2)$  of crossed modules  form $\alpha\colon A\rightarrow B$ to $\alpha'\colon A'\rightarrow B'$ in $\TC$ as {\em cover} if  $f_1\colon A\rightarrow A'$ is an isomorphism in $\TC$.

	Let $G$ be a topological internal groupoid, i.e., an object of $\Cat(\TC)$.  Let    ${\Cov}_{\Cat(\TC)}/G$  be the category of covers of $G$ in the category $\Cat(\TC)$. So the objects of ${\Cov}_{\Cat(\TC)}/G$  are  the
	covering morphisms $p\colon \widetilde{G}\rightarrow G$ over $G$ in $\Cat(\TC)$ and a morphism   from $p\colon \widetilde{G}\rightarrow G$ to $q\colon \widetilde{H}\rightarrow G$ is a
	morphism $f\colon \widetilde{G}\rightarrow \widetilde{H}$ in $\Cat(\TC)$ such that $qf=p$.
	
	The algebraic case of the following theorem was proved in \cite[Theorem 5.3]{Ak-Na-Mu-Tu}. We give the topological case of this theorem in $\TC$ as follows. The proof is obtained by Theorem \ref{Theocatequivalencess} and Theorem \cite[Theorem 5.3]{Ak-Na-Mu-Tu}.
	\begin{Theo} \label{Theocatequivalence}  Let  $G$ be an object of $\Cat(\TC)$ and    $\alpha\colon A\rightarrow B$  the  crossed module in $\TC$  corresponding to $G$ by  Theorem \ref{Theocatequivalencess}.  Let ${\Cov}_{\XMod(\TC)}/(\alpha\colon A\rightarrow B)$ be the category of covers of $\alpha\colon A\rightarrow B$ in $\TC$.  Then the categories ${\Cov}_{\Cat(\TC)}/G$ and ${\Cov}_{\XMod(\TC)}/(\alpha\colon A\rightarrow B)$ are equivalent. \end{Theo}

\end{document}